\documentclass[12pt]{amsart}
\usepackage{amsmath, amssymb, latexsym, amsthm, mathrsfs, color, mathtools}
\usepackage[top=2.5cm, bottom=2.5cm, left=2.5cm, right=2.5cm]{geometry}

\newcommand{\nc}{\newcommand}

\nc{\thusfar}{\begin{center}{\my{--- Edited thus far ---}}\end{center}}
\nc{\lei}{\le^\oo}
\nc{\card}[1]{\left|#1\right|}
\nc{\medcard}[1]{\biggl|\,#1\,\biggr|}
\nc{\smallmedcard}[1]{\bigl|\,#1\,\bigr|}
\nc{\smallcard}[1]{|\,#1\,|}
\nc{\bds}{bidirectional $\roth$-scale}
\nc{\bfP}{\mathbf{P}}
\nc{\bfQ}{\mathbf{Q}}
\nc{\bbT}{\mathbb{T}}
\nc{\bbZ}{\mathbb{Z}}
\nc{\bbN}{\mathbb{N}}
\nc{\bbC}{\mathbb{C}}
\nc{\beq}{\begin{eqnarray*}}\nc{\eeq}{\end{eqnarray*}}
\nc{\mbq}{\mb{?}}
\nc{\mb}[1]{{\mbox{\textbf{#1}}}}
\nc{\nop}{$\times$}
\nc{\fbn}{\!\!\fbox{\!\nop\!}\!\!}
\nc{\yup}{\checkmark}
\nc{\forces}{\Vdash}
\nc{\name}[1]{\dot{#1}}
\nc{\tf}{\my{FINISHED THUS FAR}}
\nc{\FU}{Fr\'echet--Urysohn}
\nc{\gs}{$\gamma$~space}
\nc{\Ga}{\Gamma}\nc{\Om}{\Omega}
\nc{\smallbinom}[2]{\begin{psmallmatrix} #1\\ #2 \end{psmallmatrix}}
\nc{\bgamma}{\smallbinom{\Om}{\Ga}}
\newcommand{\two}{\{0,1\}}
\nc{\productive}[2]{(#1,\allowbreak #2)^\x}
\nc{\prdct}[1]{(#1)^\x}
\nc{\Sel}{\mathsf{S}}
\nc{\sset}[2]{\{\,#1 : #2\,\}}
\nc{\smb}[1]{{\!\!\mb{#1}\!\!}}
\nc{\medset}[2]{{\biggl\{\,#1 : #2\,\biggr\}}}
\nc{\smallmedset}[2]{{\bigl\{\,#1 : #2\,\bigr\}}}
\nc{\set}[2]{{\left\{\,#1 : #2\,\right\}}}
\nc{\seq}[2]{{\la\, #1 : #2\,\ra}}
\nc{\eseq}[1]{#1_1, \allowbreak #1_2, \allowbreak\dotsc} 
\nc{\cube}{(\Cantor)^\bbN}
\nc{\Match}{\op{Match}}
\nc{\concat}[1]{\hat{\phantom{a}}\langle #1\rangle}
\nc{\poset}{\mathbb{P}}
\nc{\fn}[1]{{\op{Fn}(#1\times\w,2)}}
\nc{\linadd}{\op{linadd}}
\nc{\nonprod}{\non^\x}
\nc{\alephes}{{\aleph_0}}
\nc{\my}[1]{\marginpar{\textcolor{red}{***}}\textcolor{red}{#1}}
\nc{\later}[1]{{\color{green} #1}}
\nc{\BTs}[1]{{\color{green} #1 (BT)}}
\nc{\Cp}{\op{C}_\mathrm{p}}
\nc{\Bp}{\op{B}_p}
\nc{\Pa}[8]{\bibitem{#1} {#2}, \emph{#3}, {#4} \textbf{#5} ({#6}), {#7}--{#8}.}
\nc{\tPa}[5]{\bibitem{#1} {#2}, \emph{#3}, {#4}, to appear.}
\nc{\sPa}[4]{\bibitem{#1} {#2}, \emph{#3}, {#4}, submitted.}
\nc{\Bc}[9]{\bibitem{#1} {#2}, \emph{#3}, in: \textbf{#4} (#5), #6 #7, #8--#9.}
\nc{\fD}{\mathfrak{D}}
\nc{\fX}{\mathfrak{X}}
\nc{\Onbd}{\Op_{\mathrm{nbd}}} 
\nc{\Omnb}{\Om_{\mathrm{nbd}}} 
\nc{\od}{\mathfrak{od}}
\nc{\Setting}[7]{\xymatrix@R=4pt@C=7pt{#1\ar@{-}[r]&#2\ar@{-}[r]&#3\\&#4\ar@{-}[u]\\
#5\ar@{-}[uu]\ar@{-}[r] & #6\ar@{-}[u]\ar@{-}[r] & #7\ar@{-}[uu]}}
\nc{\mx}[1]{\begin{matrix}#1\end{matrix}}
\nc{\plim}{p\txt{-}\lim}
\nc{\Bgp}{{\Z^\bbN}}
\nc{\Cgp}{{{\Z_2}^\bbN}}
\nc{\Cite}[1]{\textbf{[#1]}}
\nc{\Next}[1]{{#1^+}}
\nc{\cFin}{\mathrm{cF}}
\nc{\scsp}{\text{-scale space}}
\nc{\cfn}{\text{cofinal}\ }
\nc{\Con}{\text{Concentrated}}
\nc{\Lind}{\text{Lindel\"of}\,}
\nc{\con}{\text{-Concentrated}}
\nc{\lind}{\text{-Lindel\"of}\,}
\nc{\ctbl}{\text{countably }\allowbreak}
\nc{\Hur}{\text{Hurewicz}}
\nc{\intvl}[2]{{[#1(#2),\allowbreak #1(#2\!+\!1))}}
\nc{\Bdd}{\mathbf{B}}
\nc{\Dfin}{\mathfrak{D}_\mathrm{fin}}
\nc{\grbl}{{\mbox{\textit{\tiny gp}}}}
\nc{\bbP}{\mathbb{P}}
\nc{\BOfat}{\B_{\Om_{\mathrm{fat}}}}
\nc{\Bgood}{\B_{\mathrm{good}}}
\nc{\compactN}{\cl{\mathbb{N}}}
\nc{\blocks}[2]{\op{cl}_{#2}(#1)}
\nc{\blocksplus}[2]{\op{cl}^+_{#2}(#1)}
\nc{\arx}[1]{\texttt{http://arxiv.org/math/#1}}
\nc{\bq}{\begin{quote}}
\nc{\eq}{\end{quote}}
\nc{\cl}[1]{\overline{#1}}
\nc{\Cl}[2]{\mathrm{cl}_{#1}(#2)}
\nc{\CH}{the Continuum Hypothesis}
\nc{\MA}{Martin's Axiom}
\nc{\Bfat}{\B_\mathrm{fat}}
\nc{\inv}{^{-1}}
\nc{\Cantor}{{\two^\bbN}}
\nc{\bP}{\mathbf{P}}
\nc{\bof}{\op{\fb}}
\nc{\dof}{\op{\fd}}
\nc{\bofF}{\bof(\cF)}
\nc{\sr}[3]{\underset{\mbox{#3}}{\mbox{#1}}}
\nc{\gp}{\binom{\Om}{\Ga}}
\nc{\gpsmall}{\mbox{$\gp$}}
\nc{\gig}{\gimel}
\nc{\gns}{\sone(\Om,\gig)}
\nc{\nsr}[2]{#1}
\nc{\Srg}{{\mathbb{S}}}
\nc{\Srgs}{{\mathbb{S}^*}}
\nc{\NN}{{\bbN^{\bbN}}}
\nc{\ZN}{{\Z^{\bbN}}}
\nc{\NNup}{{\bbN^{\uparrow\bbN}}}
\nc{\Pof}{\op{P}}
\nc{\PN}{{\Pof(\bbN)}}
\nc{\roth}{{[\bbN]^{\mbox{\tiny $\infty$}}}} 
\nc{\Fin}{\mathrm{Fin}}
\nc{\ici}{[\bbN]^{ \infty, \infty}}
\nc{\Inc}{{\compactN^{\uparrow\bbN}}}
\nc{\powInc}[1]{{\big(\Inc\big)^{#1}}}
\nc{\powFin}[1]{{\big(\Fin\big)^{#1}}}
\nc{\powPN}[1]{{\big(\PN\big)^{#1}}}
\nc{\NcompactN}{{\compactN^\bbN}}
\nc{\Uarrow}{\smash{\big\uparrow}}
\nc{\LE}{\preccurlyeq}
\nc{\GE}{\succcurlyeq}
\nc{\op}{\operatorname}
\nc{\im}{\op{im}}
\nc{\Span}{\op{span}}
\nc{\maxfin}{\op{maxfin}}
\nc{\ran}{\op{range}}
\nc{\iso}{\cong}
\nc{\Madd}{{\M}^\star}
\nc{\cI}{\mathcal{I}}
\nc{\cJ}{\mathcal{J}}
\nc{\scrA}{\mathscr{A}}
\nc{\scrB}{\mathscr{B}}
\nc{\scrC}{\mathscr{C}}
\nc{\scrD}{\mathscr{D}}
\nc{\scrF}{\mathscr{F}}
\nc{\scrK}{\mathscr{K}}
\nc{\A}{\D\forall}
\nc{\B}{\mathrm{B}}
\nc{\cB}{\mathcal{B}}
\nc{\cZ}{\mathcal{Z}}
\nc{\bB}{\mathbf{B}}
\nc{\BS}{\mathbf{B}(\mathcal{S})}
\nc{\BF}{\mathbf{B}(\mathcal{F})}
\nc{\BU}{\mathbf{B}(\mathcal{U})}
\nc{\cSp}{\mathcal{S}^+}
\nc{\cFp}{\mathcal{F}^+}
\nc{\cUp}{\mathcal{U}^+}

\nc{\BG}{\B_\Ga}
\nc{\BL}{\B_\Lambda}
\nc{\BT}{\B_\Tau}
\nc{\BTstar}{\B_{\Tau^*}}
\nc{\BO}{\B_\Om}
\nc{\DO}{\cD_\Om}
\nc{\KO}{\cK_\Om}
\nc{\CG}{C_\Ga}
\nc{\CL}{C_\Lambda}
\nc{\CT}{C_\Tau}
\nc{\CTstar}{C_{\Tau^*}}
\nc{\CO}{C_\Om}
\nc{\COgp}{C_{\Om^{\grbl}}}
\nc{\CLgp}{C_{\Lambda^{\grbl}}}
\nc{\BOgp}{\B_{\Om}^{\grbl}}
\nc{\BLgp}{\B_{\Lambda^{\grbl}}}
\nc{\sfC}{\mathsf{C}}
\nc{\sfD}{\mathsf{D}}
\nc{\bD}{\mathbf{D}}
\nc{\Tau}{\mathrm{T}}
\nc{\cA}{\mathcal{A}}
\nc{\cK}{\mathcal{K}}
\nc{\cD}{\mathcal{D}}
\nc{\cF}{\mathcal{F}}
\nc{\cS}{\mathcal{S}}
\nc{\cT}{\mathcal{T}}
\nc{\cG}{\mathcal{G}}
\nc{\cY}{\mathcal{Y}}
\nc{\J}{\mathcal{J}}
\nc{\cL}{\mathcal{L}}
\nc{\cM}{\mathcal{M}}
\nc{\cN}{\mathcal{N}}
\nc{\cH}{\mathcal{H}}
\nc{\cO}{\mathcal{O}}
\nc{\Op}{\mathrm{O}}
\nc{\rmA}{\mathrm{A}}
\nc{\rmF}{\mathrm{F}}
\nc{\rmB}{\mathrm{B}}
\nc{\rmD}{\mathrm{D}}
\nc{\rmP}{\mathrm{P}}
\nc{\cC}{\mathcal{C}}
\nc{\cP}{\mathcal{P}}
\nc{\bbQ}{\mathbb{Q}}
\nc{\bbR}{\mathbb{R}}
\nc{\cU}{\mathcal{U}}
\nc{\cQ}{\mathcal{Q}}
\nc{\Un}{\bigcup}
\nc{\cV}{\mathcal{V}}
\nc{\cW}{\mathcal{W}}
\nc{\Z}{{\mathbb Z}}
\nc{\Impl}{\Rightarrow}
\long\def\forget#1\forgotten{\marginpar{\textcolor{green}{Forgetting...}}}
\nc{\ft}{\mathfrak{t}}
\nc{\fb}{\mathfrak{b}}
\nc{\fc}{\mathfrak{c}}
\nc{\fd}{\mathfrak{d}}
\nc{\fg}{\mathfrak{g}}
\nc{\oo}{\infty}
\nc{\fr}{\mathfrak{r}}
\nc{\fk}{\mathfrak{k}}
\nc{\bidi}{\mathfrak{bidi}}
\nc{\fu}{\mathfrak{u}}
\nc{\fh}{\mathfrak{h}}
\nc{\fp}{\mathfrak{p}}
\nc{\fj}{\mathfrak{j}}
\nc{\fs}{\mathfrak{s}}
\nc{\w}{\omega}
\nc{\x}{\times}
\nc{\Iff}{\Leftrightarrow}
\newcommand\comp{^{\text{\tt c}}}
\nc{\nin}{\notin}
\nc{\cat}{\hat{\ }}
\nc{\sub}{\subseteq}
\nc{\spst}{\supseteq}
\nc{\sm}{\setminus}
\nc{\as}{\subseteq^*}
\nc{\les}{\le^*}
\nc{\leinf}{\le^{\infty}}
\nc{\leS}{\le_S}
\nc{\leF}{\le_F}
\nc{\leU}{\le_U}
\nc{\rest}{\restriction}

\nc{\la}{\langle}
\nc{\ra}{\rangle}
\nc{\E}{\exists}
\nc{\dom}{\op{dom}}
\nc{\cov}{\op{cov}}
\nc{\add}{\op{add}}
\nc{\addmen}{\add(\mathrm{Menger})}
\nc{\cof}{\op{cof}}
\nc{\cf}{\op{cf}}
\nc{\non}{\op{non}}
\nc{\unif}{\op{non}}
\nc{\COV}{\op{COV}}
\nc{\ADD}{\op{ADD}}
\nc{\COF}{\op{COF}}
\nc{\NON}{\op{NON}}
\nc{\impl}{\to}
\nc{\Lp}{\mathcal{L_\p}}
\nc{\Wlog}{without loss of generality}
\newtheorem{thm}{Theorem}[section]
\nc{\bthm}{\begin{thm}} \nc{\ethm}{\end{thm}}
\newtheorem{prop}[thm]{Proposition}
\nc{\bprp}{\begin{prop}} \nc{\eprp}{\end{prop}}
\newtheorem{fact}[thm]{Fact}
\nc{\bfct}{\begin{fact}} \nc{\efct}{\end{fact}}
\newtheorem{prob}[thm]{Problem}
\nc{\bprb}{\begin{prob}} \nc{\eprb}{\end{prob}}
\newtheorem{lem}[thm]{Lemma}
\nc{\blem}{\begin{lem}} \nc{\elem}{\end{lem}}
\newtheorem{app}[thm]{Application}
\nc{\bapp}{\begin{app}} \nc{\eapp}{\end{app}}
\newtheorem{claim}[thm]{Claim}
\nc{\bclm}{\begin{claim}} \nc{\eclm}{\end{claim}}
\newtheorem{cor}[thm]{Corollary}
\nc{\bcor}{\begin{cor}} \nc{\ecor}{\end{cor}}
\newtheorem{conj}[thm]{Conjecture}
\nc{\bcnj}{\begin{conj}} \nc{\ecnj}{\end{conj}}
\theoremstyle{definition}
\newtheorem{defn}[thm]{Definition}
\nc{\bdfn}{\begin{defn}} \nc{\edfn}{\end{defn}}
\newtheorem{obs}[thm]{Observation}
\nc{\bobs}{\begin{obs}} \nc{\eobs}{\end{obs}}
\theoremstyle{remark}
\newtheorem{rem}[thm]{Remark}
\nc{\brem}{\begin{rem}} \nc{\erem}{\end{rem}}
\newtheorem{cnv}[thm]{Convention}
\nc{\bcnv}{\begin{cnv}} \nc{\ecnv}{\end{cnv}}
\newtheorem{exam}[thm]{Example}
\nc{\bexm}{\begin{exam}} \nc{\eexm}{\end{exam}}
\nc{\bpf}{\begin{proof}} \nc{\epf}{\end{proof}}
\nc{\be}{\begin{enumerate}}
\nc{\ee}{\end{enumerate}}
\nc{\bi}{\begin{itemize}}
\nc{\bimy}{\my{\begin{itemize}}
\nc{\eimy}{\end{itemize}}}
\nc{\itm}{\item}
\nc{\ei}{\end{itemize}}
\nc{\Subsection}[1]{\goodbreak\subsection*{#1}}

\nc{\sone}{\mathsf{S}_1}
\nc{\sfin}{\mathsf{S}_\mathrm{fin}}
\nc{\ufin}{\mathsf{U}_\mathrm{fin}}
\nc{\Split}{\mathsf{Split}}

\nc{\gone}{\mathsf{G}_1}    \nc{\gfin}{\mathsf{G}_\mathrm{fin}}

\nc{\men}{\sfin(\cO,\cO)}
\nc{\sch}{\ufin(\cO,\Omega)}
\nc{\rothb}{\text{Rothberger}}
\nc{\pmen}{\sfin(\Omega,\Omega)}
\nc{\Rothb}{\sone(\cO,\cO)}

\nc{\prothb}{\sone(\Omega,\Omega)}
\nc{\tU}{{\tilde{U}}}
\nc{\tF}{{\tilde{F}}}
\nc{\tY}{{\tilde{Y}}}
\nc{\td}{{\tilde{d}}}
\nc{\tz}{{\tilde{z}}}
\nc{\tx}{{\tilde{x}}}
\nc{\cfd}{\cf(\fd)}

\nc{\msep}{\sfin(\cD,\cD)}
\nc{\rsep}{\sone(\cD,\cD)}
\nc{\cft}{\sfin(\Omega_{\mathbf{0}},\Omega_{\mathbf{0}})}
\nc{\scft}{\sone(\Omega_{\mathbf{0}},\Omega_{\mathbf{0}})}

\nc{\Umen}{U\text{-Menger}}
\nc{\hur}{\ufin(\cO,\Gamma)}
\nc{\tUmen}{\tU\text{-Menger}}

\nc{\Men}{\text{Menger}}
\nc{\Sch}{\text{Scheepers}}

\title{Finite powers and products of Menger sets}

\author[P.~Szewczak]{Piotr Szewczak}
\address{Piotr Szewczak, Institute of Mathematics, Faculty of Mathematics and Natural Science College of Sciences, Cardinal Stefan Wyszy\'nski University in Warsaw, W\'oycickiego 1$\slash$3, 01--938 Warsaw, Poland, and Department of Mathematics, Bar-Ilan University, Ramat Gan 5290002, Israel
}
\email{p.szewczak@wp.pl}
\urladdr{http://piotrszewczak.pl}

\author[B.~Tsaban]{Boaz Tsaban}
\address{Boaz Tsaban, Department of Mathematics, Bar-Ilan University, Ramat Gan 5290002, Israel}
\email{tsaban@math.biu.ac.il}
\urladdr{http://math.biu.ac.il/~tsaban}

\author[L.~Zdomskyy]{Lyubomyr Zdomskyy}
\address{Lyubomyr Zdomskyy, Kurt G\"odel Research Center for Mathematical Logic, University of Vienna, W\"ahringer Stra\ss e 25, A-1090 Wien, Austria.}
\email{lzdomsky@gmail.com}
\urladdr{http://www.logic.univie.ac.at/\~{}lzdomsky/}

\begin{document}
	
\begin{abstract}
We construct, using mild combinatorial hypotheses,
a real Menger set that is not Scheepers, 
and two real sets that are Menger in all finite powers, with a non-Menger product.
By a forcing-theoretic argument, we show that the same holds in the Blass--Shelah model for arbitrary values of the ultrafilter and dominating number.
\end{abstract}
	
\subjclass[2010]{
Primary: 54D20; 
Secondary: 03E17. 
}
	
\keywords{Menger property, Scheepers property, additivity number, products concentrated sets, reaping number, scales.}
	
\maketitle

\section{Introduction}

By \emph{space} we mean a topological space.
A space is \emph{$\Men$} if, for each sequence $\eseq{\cU}$ of open covers of that space, there are finite sets $\cF_1\sub\cU_1,\cF_2\sub\cU_2,\dotsc$ such that the family $\Un_{n\in\bbN}\cF_n$ covers the space.
If, in this definition, we request that each finite subset of the space is contained in a set $\Un\cF_n$, then the space is \emph{$\Sch$}.
We have the following implications between considered properties of spaces.

\[
\text{All finite powers are $\Men$} \longrightarrow \Sch \longrightarrow \Men
\]
The first of these implications follows from the results of Just, Miller, Scheepers and Szeptycki~\cite[Theorem~3.9, implications in Figure~3]{coc2}.

A \emph{set of reals} is a space homeomorphic to a subspace of the real line, with the standard topology.
We restrict our consideration to the realm of sets of reals.
The three properties are consistently equivalent~\cite{ZdoMillMProd}.
In other words, for these properties to differ, special set theoretic hypotheses are necessary.
The natural hypotheses to this end concern cardinal characteristics of the continuum~\cite{blass};
the necessary definitions are provided in the next sections.
Assuming  $\cov(\cM)=\cof(\cM)$, there is a $\Men$ set of reals that is not $\Sch$ (\cite[Theorem~32]{cbc},~\cite[Theorem~2.8]{coc2}). The methods used for that are category theoretic,
and do not lend themselves for generalizations.
We obtain the same result by a purely combinatorial approach, 
using the far milder hypothesis $\fd\leq\fr$ (Theorem~\ref{thm:MNonSch}).

Assuming $\cov(\cM)=\fc$ or $\fb=\fd$, there are sets of reals with all finite powers Menger, 
whose product is not Menger (\cite[Theorem~3.3]{TsAdd},~\cite[Proposition~3.4]{Msep}).
The assumptions $\cov(\cM)=\fc$ and $\fb=\fd$ each imply $\fd\leq\fr$.
The inequality $\fd\leq\fr$ has recently played a central role in the construction of two 
Menger sets of reals whose product is not Menger~\cite{ST}.
We refine these methods to establish the mentioned stronger result from
the same hypothesis $\fd\leq\fr$, assuming that the cardinal number $\fd$ is regular (Theorem~\ref{thm:addpmen}). 
This additional assumption follows from the earlier hypothesis $\fb=\fd$.

These results, that are optimal for the used methods, suggest the question of the necessity of the hypothesis 
$\fd\leq\fr$~\cite{ST}. 
Arguing directly in the Blass--Shelah model~\cite{BlaShe89}, we answer this question in the negative.

Finally, we apply our results to products of function spaces with the topology of pointwise convergence.

\section{The main results}

\subsection{Separation of the Menger and Scheepers properties}\label{sec:PN}

Let $\bbN$ be the set of positive natural numbers, and $\roth$ be the set of infinite subsets of $\bbN$.
We identify each set $a\in\roth$ with its increasing enumeration, an element of the Baire space $\NN$.
For each natural number $n$, by $a(n)$ denote the $n$-th smallest element of the set $a$.
We have $\roth\sub\NN$, and thus every set from $\roth$ is viewed as a function.
Let $x,y\in \roth$.
The function $x$ is \emph{dominated} by the function $y$, denoted $x\les y$, if the set $\sset{n}{y(n)<x(n)}$ is finite.
A subset of $\roth$ is \emph{dominating} if 
each function in $\roth$ is dominated by some function from this set.
Let $\fd$ be the minimal cardinality of a dominating subset of $\roth$.
For a finite subset $F$ of $\roth$, define 
$\max[F]:=\sset{\max\sset{x(n)}{x\in F}}{n\in\bbN}$, an element of $\roth$. 
A subset $X$ of $\roth$ is \emph{finitely dominating} if the set 
\[
\maxfin[X]:=\sset{\max[F]}{F\text{ is a finite subset of }X}
\]
is dominating in $\roth$.
We identify the Cantor cube $\Cantor$ with the family $\PN$ of all subsets of $\bbN$, via characteristic functions.
Since the Cantor cube is homeomorphic to Cantor's set, every subspace of $\PN$ is a set of reals.
The topologies in the set $\roth$, induced from the Cantor space $\PN$, and the Baire space $\NN$, are equivalent.
We use the following combinatorial characterizations of the Menger and Scheepers properties. 
A space is Menger (Scheepers) if and only if no continuous image of the space into $\roth$ is dominating~\cite[Proposition~3]{reclaw}(finitely dominating~\cite[Theorem~7.4(2)]{SPMProd}).

A set $r\in\roth$ \emph{reaps} a family $A\sub\roth$ if, for each set $a\in A$, both sets $a\cap r$ and $a\sm r$ are infinite. 
Let $\fr$ be the minimal cardinality of a family $A\sub\roth$ that no set $r$ reaps.
A property $\bP$ of spaces is \emph{productive} if the product space of two spaces with the property $\bP$, has the property $\bP$.

\bthm\label{thm:MNonSch}
Assume that $\fd\leq \fr$.
\be
\item There is a $\Men$ set of reals that is not  $\Sch$.
\item The $\Men$ property is not productive.
\ee
\ethm

In order to prove Theorem~\ref{thm:MNonSch} and later discussion, we introduce several notions and auxiliary results.
Fix functions $x,y\in\roth$ and a set $Z\sub\roth$.
We write $x<^\infty y$ if $y\nleq^* x$ and $Z<^\infty x$ if $z<^\infty x$ for all functions $z\in Z$.
We use the latter convention to any binary relation.
Let $\max\{x,y\}:=\sset{\max\{x(n),y(n)\}}{n\in\bbN}$, an element of $\roth$.
If $x(1)\neq 1$, then define an element $\tx\in\roth$ such that  $\tx(1):=x(1)$ and $\tx(n+1):=x(\tx(n))$, for all natural numbers $n$.
For $x\in\PN$, define $x\comp:=\bbN\sm x$, and let $\ici:=\sset{a\in\roth}{a\comp\in\roth}$.
For natural numbers $n,m$, with $n<m$, define $[n,m):=\sset{i\in\bbN}{n\leq i<m}$.
For elements $a,h\in\roth$, let 
\[
a/h:=\sset{n\in\bbN}{a\cap [h(n),h(n+1))\neq\emptyset}.
\]

\blem\label{lem:mejia}
Let $Z$ be a subset of $\roth$ with $\card{Z}<\min\{\fd,\fr\}$, and $d\in\roth$.
There are elements $x,y\in\ici$ such that $Z<^\infty x,y$, and $d\les\max\{x\comp,y\comp\}$.
\elem

\bpf
For  natural numbers $n,m$ with $n<m$, let $(n,m):=\sset{i\in\bbN}{n<i<m}$.
Since $\card{Z}<\fd$, there is a function $b\in\roth$ such that the sets 
\[
I_z:=\smallmedset{n}{\smallcard{\tz/\td\cap(b(n),b(n+1))}\geq2}
\]
are infinite for all elements $z\in Z$~\cite[Theorem~2.10]{blass}.
Since $\card{Z}<\fr$, there is a set $r\in\ici$ that reaps the family $\sset{I_z}{z\in Z}$.
Let 
\[
x:=\Un_{n\in r}[\tilde{d}(b(n)),\tilde{d}(b(n+1)+1))\text{ and }y:=\Un_{n\in r\comp}[\tilde{d}(b(n)),\tilde{d}(b(n+1)+1)).
\]

Fix an element $z\in Z$ and a natural number $n\in r\comp\cap I_{z}$.
Since $\smallcard{\tz/\td\cap(b(n),b(n+1))}\geq 2$, there is a natural number $i$ such that $\tilde{z}(i),\tilde{z}(i+1)\in [\tilde{d}(b(n)+1),\tilde{d}(b(n+1)))$.
Since $x\cap [\tilde{d}(b(n)+1),\tilde{d}(b(n+1)))=\emptyset$, we have
\[
z(\tilde{z}(i))=\tilde{z}(i+1)<\tilde{d}(b(n+1))\leq x(\tilde{z}(i)).
\] 
As the set $r\comp\cap I_z$ is infinite, we have $z<^\infty x$.
Similarly, $z<^\infty y$.

Fix a natural number $k\geq\td(b(1))$. 
There is a natural number $n$ such that $k\in[\tilde{d}(b(n)),\tilde{d}(b(n+1)))$.
Assume that $n\in r$.
Since $x\comp\cap[\tilde{d}(b(n)),\tilde{d}(b(n+1))+1)=\emptyset$, 
we have \[
d(k)<d(\tilde{d}(b(n+1)))=\tilde{d}(b(n+1)+1)\leq x\comp(k).
\]
Analogously, if $n\in r\comp$, then $d(k)<y\comp(k)$.
Thus, $d\les\max\{x\comp,y\comp\}$.
\epf

Let $\kappa$ be an uncountable cardinal number.
A subset $A$ of $\roth$ is \emph{$\kappa$-unbounded}~\cite[Definition~2.1]{ST} if $\card{A}\geq \kappa$ and, for each function $x\in\roth$, we have $\card{\sset{a\in A}{a\les x}}<\kappa$.
Let $\tau\colon \PN \to \PN$ be a homeomorphism defined by $\tau(x):=x\comp$.
Let $\Fin$ be the set of finite subsets of $\bbN$.
Then $\PN=\Fin\cup\roth$.

\begin{proof}[Proof of Theorem~\ref{thm:MNonSch}]
(1) Let $\sset{d_\alpha}{\alpha<\fd}$ be a dominating set in $\roth$.
Fix an ordinal number $\alpha<\fd$.
Let $x_\alpha, y_\alpha\in\ici$ be elements obtained from Lemma~\ref{lem:mejia}, applied to the set $\sset{d_\beta}{\beta<\alpha}$ and to the function $d_\alpha$.
Since the set $X:=\sset{x_\alpha,y_\alpha}{\alpha<\fd}$ is $\fd$-unbounded, the set  $X\cup\Fin$ is $\Men$~\cite[Corollary~2.4]{ST}. 
The continuous image $\tau[X\cup\Fin]$ of the set $X\cup\Fin$ is a subset of $\roth$.
Since $\maxfin[\tau[X\cup\Fin]]$ contains a dominating set $\sset{\max\{{x_\alpha}\comp,{y_\alpha}\comp\}}{\alpha<\fd}$, the image $\tau[X\cup\Fin]$ is finitely dominating in $\roth$.
A set of reals is $\Sch$ if no continuous image of this set into $\roth$ is finitely dominating~\cite[Theorem~2.1]{dhm}.
Thus, the set $X\cup\Fin$ is not $\Sch$.

(2) Since the continuous image of the product space $(X\cup\Fin)^2$ under the map $(x,y)\mapsto\max\{x\comp,y\comp\}$ is dominating in $\roth$, the product space $(X\cup\Fin)^2$ is not $\Men$.
\epf

\brem
Let $\cov(\cM)$ be the minimal cardinality of a family of meager subsets of $\roth$ that covers $\roth$ and $\cof(\cM)$ be the minimal cardinality of a \emph{cofinal} family of meager sets in $\roth$, i.e., every meager subset of $\roth$ is contained in a member of the family.
Assuming $\cov(\cM)=\cof(\cM)$, there is a $\Men$ set  of reals that is not $\Sch$~\cite[Theorem~32]{cbc}.
Since $\cov(\cM)\leq\fd\leq\cof(\cM)$ and $\cov(\cM)\leq \fr$, the equality $\cov(\cM)=\cof(\cM)$ implies that $\fd\leq\fr$.
Thus, Theorem~\ref{thm:MNonSch}(1) is a substantial extension of the mentioned result.
\erem

\brem\label{rem:MenNonSch}
It is always the case that, if there is a $\Men$ set of reals that is not $\Sch$, then the $\Men$ property is not productive:
Let $X$ be such a set.
Some finite power of $X$ is not $\Men$~\cite[Theorem~3.9]{coc2} and let $n$ be the minimal natural number with this property. 
The product space of $\Men$ sets $X$ and $X^{n-1}$ is not $\Men$.
\erem

\subsection{Products of $\Men$ sets with strong properties}

A property of spaces $\bP$ is \emph{additive} if the union of any two spaces with the property $\bP$, has the property $\bP$.

\bthm \label{thm:addpmen}
Assume that $\fd\leq \fr$ and the cardinal number $\fd$ is regular.
\be
\item There are two sets of reals whose all finite powers are $\Men$, but whose union is not $\Sch$.
\item The $\Sch$ property is not additive.
\item There are two sets of reals whose all finite powers are $\Men$, but whose product space is not $\Men$.
\item The $\Sch$ property is not productive.
\ee
\ethm

For elements $x,y\in\roth$, let $[x\leq y]:=\sset{n}{x(n)\leq y(n)}$ and $[x< y]:=\sset{n}{x(n)< y(n)}$.
Let $A$ be a subset of $\roth$, and $x,y\in\roth$.
We write $x\leq_A y$ if $[x\leq y]\in A$, and $x<_A y$ if $[x< y]\in A$.
A subset $X$ of $\roth$ is \emph{$\le_A$-bounded} if there is a function $b\in\roth$ such that $X\le_A b$. A subset of $\roth$ is \emph{$\le_A$-unbounded} if it is not $\le_A$-bounded.

A  \emph{filter} is a subset of $\roth$ with empty intersection, that is closed under finite intersections and taking supersets.
An \emph{ultrafilter} is a maximal filter.
Let $U$ be an ultrafilter, and $\bof(U)$ be the minimal cardinality of a $\leU$-unbounded set in $\roth$.
A set $X\sub \roth$ is a \emph{$U$-scale}~\cite[Definition~4.1]{ST} if $\card{X}\geq \bof(U)$ and, for each function $b\in\roth$, we have $\card{\sset{x\in X}{x\leU b}}<\bof(U)$.
A space $X$ is \emph{$\Umen$} if, for each sequence $\eseq{\cU}$ of open covers of the space $X$, there are finite sets $\cF_1\sub\cU_1,\cF_2\sub\cU_2,\dotsc$ such that, for each point $x\in X$, we have $\sset{n}{x\in\Un\cF_n}\in \cU$.
A set of reals is $\Umen$ if every continuous image of the set into $\roth$ is $\leU$-bounded~\cite[Proposition~4.4]{ST}.
Since no $\leU$-bounded set in $\roth$ is finitely dominating, every $\Umen$ set is $\Sch$~\cite[Theorem~7.4(2)]{SPMProd}.

An ultrafilter $U$ is \emph{near coherent} to an ultrafilter $\tU$ if there is a function $h\in\roth$ such that for each set $u\in U$, we have 
\[
\Cl{h}{u}:=\Un\smallmedset{[h(n),h(n+1))}{u\cap [h(n),h(n+1))\neq\emptyset}\in\tU.
\]
The near coherence relation is an equivalence relation.
The near coherence of filters (NCF) is the statement introduced by Blass that any two ultrafilters are near coherent~\cite{blassNCF}.
NCF is independent of ZFC~\cite{blassNCF} and it follows, e.g.,  from the inequality $\fu<\fg$~\cite{laflamme}, where $\fu$ is the ultrafilter number and $\fg$ is the groupwise density number.
It demonstrates that the forthcoming Lemma~\ref{lem:NonCoh} needs a set theoretic assumption.

In order to prove Theorem~\ref{thm:addpmen}, we need the following Lemmata.
Let $\cf(\fd)$ be the cofinality of the cardinal number $\fd$.
A subset of $\roth$ is \emph{centered} if the finite intersections of its elements, are infinite.
For  natural numbers $n,m$ with $n\leq m$, let $[n,m]:=\sset{i\in\bbN}{n\leq i\leq m}$.
For elements $x,y\in\roth$, we write $x\leq y$ if $[x\leq y]=\bbN$. 

\blem \label{lem:NonCoh}
Assume that $\fd\leq\fr$. 
There are non near coherent ultrafilters $U$ and $\tU$ such that $\bof(U)=\bof(\tU)=\cf(\fd)$. 
\elem

\bpf
Let $\seq{D_\alpha}{\alpha<\cfd}$ be a strictly increasing sequence of subsets of $\roth$ of cardinality less than $\fd$ such that the set  $\Un_{\alpha<\cfd}D_\alpha$ is dominating in $\roth$.
We may assume that for each function $x\in \roth$ there is a function $d$ from the dominating set with $x\leq d$.
Proceed by transfinite induction on ordinal numbers $\alpha<\fd$.
Define increasing sequences $\seq{F_\alpha}{\alpha<\fd}$, $\seq{\tF_\alpha}{\alpha<\fd}$ of subsets of $\roth$, and a set $\sset{x_\alpha}{\alpha<\cfd}\sub\ici$ such that, for each ordinal number $\alpha<\fd$:
\bi
\item the sets $F_\alpha$, $\tF_\alpha$ with $\smallcard{F_\alpha},\smallcard{\tF_\alpha}<\fd$ are closed under finite intersections,
\item $D_\alpha\cup\sset{x_\beta}{\beta<\alpha}\le_{F_\alpha} x_\alpha$ and $D_\alpha\cup\sset{{x_\beta}\comp}{\beta<\alpha}\le_{\tF_\alpha} {x_\alpha}\comp$.
\ei

For functions $s,f\in\roth$, let $s\circ f$ be a function in $\roth$ such that $(s\circ f)(n):=s(f(n))$, for all natural numbers $n$.
Fix an ordinal number $\alpha<\fd$.
Let
\begin{align*}
Y_\alpha:=&\sset{s \circ f}{s\in \maxfin[D_\alpha\cup\sset{x_\beta}{\beta<\alpha}], f\in\Un_{\beta<\alpha}F_\beta},\\
\tY_\alpha:=&\sset{s \circ f}{s\in \maxfin[D_\alpha\cup\sset{{x_\beta}\comp}{\beta<\alpha}], f\in\Un_{\beta<\alpha}\tF_\beta}.
\end{align*}

Since $\smallcard{Y_\alpha\cup\tY_\alpha}<\fd\leq\fr$, there is an element $x_\alpha\in\ici$ such that $Y_\alpha\cup \tY_\alpha<^\infty x_\alpha$ and $Y_\alpha\cup \tY_\alpha<^\infty {x_\alpha}\comp$~\cite[Lemma~3.4]{ST}.

Fix elements $s\in \maxfin[D_\alpha\cup\sset{x_\beta}{\beta<\alpha}]$ and $f\in\Un_{\beta<\alpha}F_\beta$. 
Since the set $[(s\circ f)\leq x_\alpha]$ is infinite, the intersection $f\cap [s\leq x_\alpha]$ is infinite, too.
Thus, the set
\[
\Un_{\beta<\alpha}F_\beta\cup\set{[s\leq x_\alpha]}{s\in \maxfin[D_\alpha\cup\sset{x_\beta}{\beta<\alpha}]}
\]
is centered. Let $F_\alpha$ be the latter set closed under finite intersections.
Then $\card{F_\alpha}<\fd$.
Similarly, define a set $\tF_\alpha$. 

The sets $F:=\Un_{\alpha<\cfd}F_\alpha$ and $\tF:=\Un_{\alpha<\cfd}\tF_\alpha$ are closed under finite intersections. 
Let $U$ and $\tU$ be ultrafilters containing the sets $F$ and $\tF$, respectively.
We have $\bof(U)=\cfd$:
Let $Y$ be a subset of $\roth$ with $\card{Y}<\cfd$. 
There is an ordinal number $\alpha<\cfd$ such that, for each function $y\in Y$, there is a function $d_y\in D_\alpha$ with $y\les d_y$. 
Since $D_\alpha\leU x_\alpha$, we have $Y\leU x_\alpha$.

The ultrafilters $U$ and $\tilde{U}$ are not near coherent:
Fix $h\in \roth$.
There are an ordinal number $\alpha<\mathfrak d$ and an element  $d\in
D_\alpha$ such that $2h(n+1)<d(h(n))$ for all natural numbers $n$.
We may assume that  $\bbN\in
F_0\cap\tilde{F}_0$.
Since $d\in D_\alpha\sub
Y_\alpha$, we have $[d\leq x_\alpha]\in F_\alpha$ and
$[d\leq x_\alpha\comp]\in \tilde{F}_\alpha$.
Next, we show that $\Cl{h}{[d\leq x_\alpha]}\cap \Cl{h}{[d\leq x_\alpha\comp]}=\emptyset$.
Suppose contrary that there is a natural number $n$ such that $[h(n), h(n+1))\cap [d\leq x_\alpha]\neq\emptyset$ and $[h(n), h(n+1))\cap [d\leq x_\alpha\comp]\neq\emptyset$. 
Pick 
\[
l\in [h(n), h(n+1))\cap [d\leq x_\alpha]\text{ and }k\in [h(n), h(n+1))\cap [d\leq x_\alpha\comp].
\]
We have
\[
 2h(n+1)<  d(h(n))\leq d(l)\leq x_\alpha(l)<x_{\alpha}(h(n+1))
\]
and
\[
 2h(n+1)<  d(h(n))\leq d(k)\leq x_\alpha\comp(k)<x_{\alpha}\comp(h(n+1)).
\]
From the above equations, we get
\[
 [1,2h(n+1)+1]\sub
 x_\alpha([1,h(n+1)])\cup
 x_\alpha\comp([1,h(n+1)]),
\]
and thus
\begin{gather*}
2h(n+1)+1= \smallmedcard{[1,2h(n+1)+1]}\leq 
 \smallmedcard{x_\alpha([1,h(n+1)])}+\smallmedcard{
 x_\alpha\comp([1,h(n+1)])}=\\
 h(n+1)+h(n+1)=2h(n+1),
\end{gather*}
a contradiction.
\epf

\brem
The non near coherence of the ultrafilters $U$ and $\tilde{U}$ in the proof of Lemma~\ref{lem:NonCoh} could be also established with the help of the following shorter topological
argument, which, however, is not self-contained.

Let $\sset{x_\alpha}{\alpha<\cf(\fd)}$ be defined as in the proof of Lemma~\ref{lem:NonCoh}.
The set $\sset{x_\alpha}{\alpha<\cf(\fd)}$ is a $U$-scale:
Let $b\in\roth$. There is an ordinal number $\alpha<\cfd$ such that the function $b$ is dominated by some function from the set $D_\alpha$. 
For each ordinal number $\beta$ with $\alpha\leq\beta<\cfd$, we have $D_\alpha\leU x_\beta$, and thus $b\leU x_\beta$.
Analogously, we have $\bof(\tU)=\cfd$ and the set $\sset{{x_\alpha}\comp}{\alpha<\cf(\fd)}$ is a $\tU$-scale.
Let $X:=\sset{x_\alpha}{\alpha<\fd}$.
Then the set $X\cup\Fin$ is $\Umen$~\cite[Theorem~5.3(2)]{ST}. The continuous image $\tau[X\cup\Fin]$ of the set $X\cup\Fin$ in $\roth$ contains a $\le_{\tU}$-unbounded set $\sset{{x_\alpha}\comp}{\alpha<\fd}$, and thus the set $X\cup\Fin$ is not $\tUmen$~\cite[Proposition~4.4]{ST}. 
We conclude that the ultrafilter $U$ is not near coherent to  the ultrafilter $\tU$~\cite[Theorem~2.32]{sfh}.
\erem

\blem\label{lem:NonCohGap}
Let $U$ and $\tU$ be non near coherent ultrafilters, and $d\in\roth$. There are sets $I,J\in\roth$ such that 
\[
\Un_{n\in I}[d(n),d(n+1))\in U, \Un_{n\in J}[d(n),d(n+1))\in \tU,\text{ and }\card{i-j}>3, \text{ for all }i\in I, j\in J.
\]
\elem

\bpf
Let $b_k:=\sset{8n+k}{n\in\bbN}$, where $k\in\{0,\dotsc,7\}$.
Let $k_0$ be such that $x:=\bigcup_{i\in b_{k_0}}[d(i), d(i+1))\in U$.
Define $d':=\sset{d(n)}{|n-i|>3 \text{ for all }i\in b_{k_0}}$, an element of $\roth$.
Since the ultrafilters $U$ and $\tU$ are non near coherent, there are sets $y\in U$ and $\tilde{y}\in\tU$
such that $\Cl{d'}{y}\cap \Cl{d'}{\tilde{y}}=\emptyset$.
We may assume that $y\sub x$.
Let $I:=\sset{n}{[d(n),d(n+1))\cap y\neq\emptyset}$
and
$J:=\sset{n}{[d(n),d(n+1))\cap \tilde{y}\neq\emptyset}$.
Then $\Un_{n\in I}[d(n),d(n+1))\in U$ and $\Un_{n\in J}[d(n),d(n+1))\in \tU$.
Fix $i\in I$ and $j\in J$.
Since $y\sub x$, we have $i\in b_{k_0}$.
Let $n, m$ be natural numbers such that $d(n),d(m)$ are consecutive elements of $d'$ with $n\leq j<j+1\leq m$.
By the definition of $d'$, we have $|n-i|, |m-i|>3$.
Since $[d(n),d(m))\cap y=\emptyset$, we have $i<n$ or $m\leq i$, and thus $|i-j|>3$.
\epf

\blem \label{lem:NonCohStep}
Let $U$ and $\tU$ be non near coherent ultrafilters, and $d\in\roth$.
There are elements $x,y\in \ici$ such that $d<_U x$,
$d<_\tU y$, and $d\leq^*\max\{x\comp,y\comp\}$.
\elem

\bpf
By Lemma~\ref{lem:NonCohGap}, there are sets $I,J\in\roth$ such that 
\[
\Un_{n\in I}[\td(n),\td(n+1))\in U, \Un_{n\in J}[\td(n),\td(n+1))\in \tU,\text{ and }\card{i-j}>3
\]
for all natural numbers $i\in I$, and $j\in J$.
Let $x,y\in\roth$ be elements such that 
\[
x\comp:=\Un_{n\in I}[\td(n),\td(n+2)),\text{ and }y\comp:=\Un_{n\in J}[\td(n),\td(n+2)).
\]

Fix a natural number $n\in I$ and a natural number $k\in[\td(n),\td(n+1))$.
Since $x\cap[\td(n),\td(n+2))=\emptyset$, we have $d(k)< d(\td(n+1))=\td(n+2)\leq x(k)$.
Thus, $d<_U x$. 
Similarly, we have $d<_\tU y$. 

Fix a natural number $n$ and a natural number $k\in[\td(n),\td(n+1))$.
One of the sets $x\comp$, or $y\comp$ has empty intersection with the interval $[\td(n),\td(n+2))$.
If $x\comp\cap [\td(n),\td(n+2))=\emptyset$, then 
$d(k)<d(\td(n+1))<\td(n+2)\leq x\comp(k)$.
If the set $y\comp$ has this property, then $d(k)\leq y\comp(k)$.
Thus, $d\les\max\{x\comp,y\comp\}$
\epf


Let $\bP$ be a property of spaces.
A space $X$ is \emph{productively $\bP$} if, for each space $Y$ with the property $\bP$, the product space $X\x Y$, has the property $\bP$.

\blem\label{lem:Uscales}
Assume that $\fd\leq\fr$. 
Let $U$ be an ultrafilter with $\bof(U)=\fd$.
For any ultrafilter $\tU$ non near coherent to the ultrafilter $U$ with $\bof(\tU)=\fd$, there are a $U$-scale $X\sub\roth$, and a $\tU$-scale $Y\sub\roth$ such that
\be
\item the union $(X\cup\Fin)\cup(Y\cup\Fin)$ is not $\Sch$,
\item the product space $(X\cup\Fin)\x (Y\cup\Fin)$ is not $\Men$.
\ee
\elem

\bpf
(1)
Let $\sset{e_\alpha}{\alpha<\fd}$ be a dominating set in $\roth$.
Proceed by transfinite induction on ordinal numbers $\alpha<\fd$.
Let $d_0=e_0$.
Fix an ordinal number $\alpha<\fd$ and assume that elements $d_\beta\in\roth$ have been already defined for all ordinal numbers $\beta$  with $\beta<\alpha$. Since $\bof(U)=\bof(\tU)=\fd$, there is a function  $d_\alpha\in\roth$ such that $\sset{d_\beta}{\beta<\alpha}\leU d_\alpha$, $\sset{d_\beta}{\beta<\alpha}\le_{\tU} d_\alpha$ and $e_\alpha\les d_\alpha$.
Then the set $\sset{d_\alpha}{\alpha<\fd}$ is a dominating set in $\roth$, simultaneously a $U$-scale, and a $\tU$-scale.

For each ordinal number $\alpha<\fd$, let $x_\alpha,y_\alpha\in\roth$ be elements obtained from Lemma~\ref{lem:NonCohStep}, applied to the function $d_\alpha$.
The sets $X:=\sset{x_\alpha}{\alpha<\fd}$, and  $Y:=\sset{y_\alpha}{\alpha<\fd}$ are a $U$-scale, and a $\tU$-scale in $\roth$, respectively.
Thus the sets $X\cup\Fin$, and $Y\cup\Fin$ are productively $\Umen$, and productively $\tUmen$, respectively~\cite[Theorem~5.3(2)]{ST}.
Let $Z:=X\cup Y\cup\Fin$.
The set $\maxfin[\tau[Z]]$, a subset of $\roth$, contains a dominating set $\sset{\max\{{x_\alpha}\comp,{y_\alpha}\comp\}}{\alpha<\fd}$.
The continuous image $\tau[Z]$ of the set $Z$ is finitely dominating in $\roth$, and thus the set $Z$ is not $\Sch$~\cite[Theorem~2.1]{dhm}. 

(2) The continuous image of the product space $(X\cup\Fin)\x (Y\cup\Fin)$, under the map $(x,y)\mapsto \max\{x\comp,y\comp\}$, is dominating in $\roth$. 
Thus, this product space is not $\Men$.
\epf

\begin{proof}[Proof of Theorem~\ref{thm:addpmen}]

(1) Since  the cardinal number $\fd$ is regular, there is an ultrafilter $U$ with $\bof(U)=\fd$~\cite{canjar}. By Lemma~\ref{lem:NonCoh} and the fact that the near  coherence relation is an equivalence relation, there is an ultrafilter $\tU$ non near coherent to the ultrafilter $U$ with $\bof(\tU)=\fd$.
Let $X$ and $Y$ be a $U$-scale, and a $\tU$-scale, respectively, obtained from Lemma~\ref{lem:Uscales}.
The sets $X\cup\Fin$, and $Y\cup\Fin$, are productively $\Umen$, and productively $\tUmen$, respectively, and thus all finite powers of these sets are $\Men$. Apply Lemma~\ref{lem:Uscales}(1).

(2) If all finite powers of a space are $\Men$, then the space is $\Sch$. Apply (1).

(3) Apply sets from (1) and  Lemma~\ref{lem:Uscales}(2).	

(4) If all finite powers of a space are $\Men$, then the space is $\Sch$. Apply (3).
\epf

It is consistent with ZFC, that for any ultrafilter $U$, the $\Umen$ and $\Sch$ properties are equivalent~\cite[proof of Theorem~3.7]{semtrich}.
In contrast to this result, we have the following corollary from Theorem~\ref{thm:addpmen}.

\bcor\label{cor:addpmen}
Assume that $\fd\leq\fr$. 
For each ultrafilter $U$, there is a $\Sch$ set of reals that is not $\Umen$.
\ecor

\bpf
Assume that $\bof(U)<\fd$.
Since every set of reals of cardinality less than $\fd$ is $\Sch$~\cite[Theorem~2.1]{dhm}, each $\leU$-unbounded subset of $\roth$ of cardinality $\bof(U)$, is $\Sch$, but not $\Umen$.
Assume that $\bof(U)=\fd$.
By Lemma~\ref{lem:Uscales}(2), there are an ultrafilter $\tU$, and two sets of reals that are productively $\Umen$, and productively $\tUmen$, respectively, whose product space is not $\Men$. 
One of these sets is $\tUmen$, and thus $\Sch$~\cite[Proposition~4.4]{ST},~\cite[Theorem~2.1]{dhm}, but it is not $\Umen$.
\epf

\subsection{Additivity and productivity of the $\Men$ property.}
Let $\addmen$ be the minimal cardinality of a family of $\Men$ sets of reals, whose union is not $\Men$.
For every subset $X$ of $\roth$, containing a $\cfd$-unbounded, or a  $\fd$-unbounded set, there is a $\Men$ set of reals $Y$ such the product space $X\x Y$ is not $\Men$~\cite[Theorem~2.7]{ST}.
We prove the same statement, assuming that the set $X$ contains an $\addmen$-unbounded set.

\bthm\label{thm:addMunbdd}
Let $X$ be a subset of $\roth$ containing an $\addmen$-unbounded set.
There is a $\Men$ set of reals $Y$ such that the product space $X\x Y$ is not $\Men$.
\ethm

Recall the result from an earlier work~\cite{ST}, a main tool in the proof of Theorem~\ref{thm:addMunbdd}.

\blem[{\cite[Lemma~2.6]{ST}}]\label{lem:ergosum}
For sets $a,b\in\PN$, let
\[
a\uplus b:=(2a)\cup(2b+1)=\set{2k}{k\in a}\cup\set{2k+1}{k\in b}.
\]
Then:
\be
\item For each set $a\in\roth$ and each natural number $n$, we have 
$(a\uplus a)(2n) =2a(n)+1$.
\item For all sets $a,b,c,d\in\roth$ with $a\le b$ and $c\le d$, we have $a\uplus c\le b\uplus d$.\qed
\ee
\elem

For sets $x,y$, let $x\oplus y:=(x\cup y)\sm(x\cap y)$.
The space $\PN$, with the group
operation $\oplus$ and the standard topology (see the beginning of
Subsection~\ref{sec:PN} for the explanation of what kind of topology on
$\PN$ we consider) is a topological group.
For sets $X,Y\sub\PN$, let $X\oplus Y:=\sset{x\oplus y}{x\in X, y\in Y}$.
For an element $x\in \roth$ let $2x:=\sset{2k}{k\in x}$, and $(x+1):=\sset{k+1}{k\in x}$.
For a set $X\sub\roth$ and an element $a\in\roth$, let $2X:=\sset{2x}{x\in X}$, $X+1:=\sset{x+1}{x\in X}$ and $a\uplus X:=\sset{a\uplus x}{x\in X}$.

\begin{proof}[Proof of Theorem~\ref{thm:addMunbdd}]
Let $\kappa:=\addmen$.
There is a family $\sset{M_\alpha}{\alpha<\kappa}$ of Menger sets of reals whose union is not Menger.
Then there is a continuous image of this union in $\roth$, a dominating set.
Since the Menger property is preserved by continuous mappings, we may assume that $\Un_{\alpha<\kappa}M_\alpha$ is a union of $\Men$ subsets of $\roth$.
Let $\sset{a_\alpha}{\alpha<\kappa}$ be a $\kappa$-unbounded subset of $X$.
For each ordinal number $\alpha<\kappa$, the set $\tilde{M}_\alpha:=\sset{\max\{a_\alpha,x\}}{x\in M_\alpha}$, a continuous image of the $\Men$ set $M_\alpha$, is $\Men$.

The set 
\[
Y:=\Un_{\alpha<\kappa}(a_\alpha\uplus \tilde{M}_\alpha)\cup\Fin
\]
is $\Men$:
Let $\eseq{\cU}$ be a sequence of open covers of the set $Y$ which are open also in $\PN$.
There are sets $U_1\in\cU_1,U_2\in\cU_2,\dotsc$ such that the set $U:=\Un_nU_n$ contains the set $\Fin$.
Since $\PN\sm U$ is a compact subset of $\roth$, there is an element $b\in\roth$ such that $\PN\sm U\leq b$.
Define $b'(n):=b(2n)$ for all natural numbers $n$.
Fix an ordinal number $\alpha<\kappa$, and an element $y\in \tilde{M}_\alpha$.
Assume that $a_\alpha\uplus y\leq b$.
By Lemma~\ref{lem:ergosum}, we have
\[
a_\alpha(n)<2a_\alpha(n)+1\leq(a_\alpha\uplus a_\alpha)(2n)\leq(a_\alpha\uplus y)(2n)\leq b(2n)=b'(n)
\]
for all natural numbers $n$.
Thus, $a_\alpha\leq b'$.
Since the set $\sset{a_\alpha}{\alpha<\kappa}$ is $\kappa$-unbounded, there is an ordinal number $\alpha'<\kappa$ such that $\sset{\beta<\kappa}{a_\beta\leq b'}\sub \alpha'$.
We conclude that the set 
\[
Y':=Y\sm U\sub\bigl(\Un_{\alpha<\kappa}a_\alpha\uplus\tilde{M}_\alpha\bigr)\cap\sset{x\in\roth}{x\leq b}
\] 
is a closed subset of the $\Men$ set $\Un_{\beta<\alpha'}(a_\beta\uplus\tilde{M}_\beta)$, and thus it is $\Men$.
There are finite sets $\cF_1\sub\cU_1, \cF_2\sub \cU_2,\dotsc$ such that the family $\Un_n\cF_n$ covers the set $Y'$.
Thus, the family $\Un_n(\cF_n\cup\{U_n\})$ is a cover of the set $Y$.

The set $2X\oplus Y$ is dominating in $\roth$:
By the definition of the operation $\uplus$, the set $2X\oplus Y$ is a subset of $\roth$.
Fix an ordinal number $\alpha<\kappa$, and an element $y\in \tilde{M}_\alpha$.
We have 
\[
2a_\alpha\oplus(a_\alpha\uplus y)=2y+1\in (2X)\oplus Y.
\]
Thus, the set $(2X)\oplus Y$ contains a dominating set $(2\Un_{\alpha<\kappa}\tilde{M}_\alpha)+1$.

Since the set $(2X)\oplus Y$ is a continuous image of the product space $X\x Y$, the product space $X\x Y$ is not $\Men$.
\epf

\bthm\label{thm:add<d}
Assume that $\addmen<\fd$.
The $\Men$ property is not productive.
\ethm

We need the following Lemma.

\blem\label{lem:addMunbdd}
There is an $\addmen$-unbounded subset of $\roth$.
\elem

\bpf
Let $\kappa:=\addmen$, and $\seq{M_\alpha}{\alpha<\kappa}$ be an increasing sequence of $\Men$ subsets of $\roth$ whose union is a dominating set.
Fix an ordinal number $\alpha<\kappa$.
The set $\sset{a_\beta}{\beta<\alpha}\cup M_\alpha$ is a union of two $\Men$ sets, and thus it is not dominating in $\roth$.
There is a function $a_\alpha\in\roth$ such that 
\[
\sset{a_\beta}{\beta<\alpha}\cup M_\alpha <^\infty a_\alpha.
\]
Let $A:=\sset{a_\alpha}{\alpha<\kappa}$, and $b\in\roth$. 
There is an ordinal number $\alpha<\kappa$, and a function $d\in M_\alpha$ such that $b\les d$.
By the construction, we have $\sset{\beta<\kappa}{a_\beta\les b}\sub\alpha$, and thus $\card{\sset{a\in A}{a\les b}}<\kappa.$ 
\epf

\begin{proof}[Proof of Theorem~\ref{thm:add<d}]
The $\addmen$-unbounded set, constructed in Lemma~\ref{lem:addMunbdd} has cardinality $\addmen$, and thus it is $\Men$. Apply Theorem~\ref{thm:addMunbdd}.
\epf

\subsection{Products of $\Men$ sets in the Blass--Shelah model}
By Theorem~\ref{thm:MNonSch}(2), assuming $\fd\leq\fr$, the $\Men$ property is not productive.
On the other hand, in the Miller model, where $\fd>\fr$, the $\Men$ property is productive~\cite{ZdoMillMProd}.
We show that the inequality $\fd\leq\fr$ is not necessary to prove that the $\Men$ property is not productive.
In this section, we shall work in the model
constructed by Blass and Shelah~\cite{BlaShe89} and use their notations.
Let $\nu$ and $\delta$ be uncountable regular ordinal numbers in a model of ZFC $+$ GCH with $\nu<\delta$.
Let $V(\delta,0)$ be the initial Cohen extension of the model of ZFC $+$ GCH, after adding $\delta$-many Cohen reals.
For each ordinal number $\xi<\nu$, let $V(\delta,\xi)$ be the model obtained after $\xi$ stages in the Mathias forcing iteration with finite support, and $s_\xi$ be a Mathias real over the model $V(\delta,\xi)$ with respect to a certain ultrafilter $U_\xi$, in this model.
Since $U_\xi$ is chosen to contain the set $\sset{s_\zeta}{\zeta<\xi}$ for all ordinal numbers $\zeta<\nu$, the sequence $\seq{s_\xi}{\xi<\nu}$ is almost decreasing.
Let $\fu$ be the minimal cardinality of a basis of an ultrafilter.
In the model $V(\delta,\nu)$, the final extension of the model of ZFC $+$ GCH, we have $\fr=\fu=\nu$, $\fd=\delta$, and thus $\fr<\fd$~\cite{BlaShe89}.

\bthm\label{thm:BSMnonProd}
In the Blass--Shelah model $V(\delta,\nu)$, there are two sets of reals whose all finite powers are $\Men$, but their product space is not $\Men$.
\ethm

Let $\fb$ be the minimal cardinality of a $\les$-unbounded set in $\roth$.
For each ordinal number $\xi<\nu$, let $M_\xi:=V(\delta,\xi)\cap\roth$.

\blem\label{lem:BSaddMen}
In the Blass--Shelah model, we have $\addmen=\fr$.
In particular, $\addmen<\fd$.
\elem

\bpf
For each ordinal number $\xi<\nu$, the set $M_\xi$ is $\Sch$ 
in the model $V(\delta,\nu)$\footnote{This is not a direct application of~\cite[Theorem~11]{SchTal10}, but rather of its proof which gives: If $X$ is a ground model set of reals and $P$ is an iteration such that each real is contained in some intermediate model, and there is an unbounded real over every intermediate model, then $X$ is Menger in the final model. }~\cite[Theorem~11]{SchTal10}. 
Thus, $\addmen\leq\nu=\fu$.
Since $\fr=\fb$~\cite[Lemma~3.3]{Mil07} and $\fb\leq\addmen$~\cite[Corollary~2.2(2)]{TsAdd}, we have $\addmen=\fr$.
\epf

\brem
Theorem~\ref{thm:add<d} and Lemma~\ref{lem:BSaddMen} imply that, in the Blass--Shelah model, the $\Men$ property is not productive.
\erem

Let $X$ be a space.
A map $\Phi\colon X \to \roth$ is \emph{upper continuous} if the sets $\sset{x\in X}{\Phi(x)(n)\leq m}$ are open for all natural numbers $n,m$.

\blem\label{lem:upper}
Let $X$ be a space and $\Phi\colon X\to \roth$ be an upper continuous map.
The map $\Phi$ is Borel.
\elem

\bpf
Let $n,m\in\bbN$.
The set
\[
\sset{x\in X}{\Phi(x)(n)= m}=\sset{x\in X}{\Phi(x)(n)\leq m}\sm \sset{x\in X}{\Phi(x)(n)< m},
\] is a difference of two open sets, and thus is Borel.
\epf

\begin{proof}[Proof of Theorem~\ref{thm:BSMnonProd}]

The family $\sset{[x<s_\xi]}{x\in M_\xi,\xi< \nu}$ is centered:
Fix a natural number $n$.
Let $\seq{\xi_i}{i\leq n}$ be a nondecreasing sequence of ordinal numbers smaller than $\nu$, and $x_i\in M_{\xi_i}$ for natural numbers $i\leq n$. 
For each ordinal number $\xi<\nu$, we have $\maxfin[M_{\xi}]=M_{\xi}$.
Thus, assume that the sequence $\seq{\xi_i}{i\leq n}$ is strictly increasing.
For each ordinal number $\xi<\nu$, and functions $x,y\in M_\xi$, by the genericity of the function $s_\xi$, the set $[x<s_\xi]\cap y$ is infinite.
Thus, the intersection $\bigcap_{i\leq n}[x_i<s_{\xi_i}]$ is infinite.

Let $U$ be an ultrafilter containing the family $\sset{[x<s_\xi]}{x\in M_\xi,\xi< \nu}$.
For each ordinal number
$\xi<\nu,$ let $\tilde{M}_\xi$ be the set $\sset{z_\xi(x)}{x\in
M_\xi}\sub\roth$, where $z_\xi(x)(n)=s_\xi(n)+1$ if
$n\in x$ and $z_\xi(x)(n)=s_\xi(n)$ otherwise. Thus $\tilde{M}_\xi$
is a continuous image of $M_\xi$ under the map $x\mapsto z_\xi(x)$,
and hence $\tilde{M}_\xi\in V_{\delta,\xi+1}$, which is going to be
used later.
Let $X:=\Un_{\xi<\nu}\tilde{M_\xi}$.

All finite powers of the set $X\cup\Fin$ are $\Men$, in the model $V(\delta,\nu)$:
We prove a formally stronger statement that, for each natural number $n$, and each finite product $M$ of the sets $M_{\xi}$, the product space $(X\cup\Fin)^n\x M$ is $\Men$.
Each finite product of the sets $M_\xi$ is $\Men$~\cite[Theorem~11]{SchTal10}.
Proceed by induction.
Fix a natural number $n$.
Let $k$ be a natural number, and $\xi_i<\nu$ be ordinal numbers for natural numbers $i\leq k$.
Let $M:=\prod_{i\leq k}M_{\xi_i}$, and $\Psi \colon (X\cup\Fin)^{n+1}\x M\to\roth$ be a continuous map.
There is an upper continuous map $\Phi\colon M \to\roth$ such that, for elements
$x=(x_1,\dotsc,x_{n+1})\in (X\cup\Fin)^{n+1}$, and $y=(y_1,\dotsc, y_k)\in M$, and a natural number $l$, we have:

\begin{center}
If $\Phi(y)(l)<\min\sset{x_j(l)}{j\leq n+1}$, then $\Psi(x,y)(l)< \Phi(y)(l)$~\cite[Lemma~5.1]{ST}.
\end{center}
By Lemma~\ref{lem:upper}, the map $\Phi$ is Borel.

Let $\xi<\nu$ be an ordinal number with $\xi>\max\sset{\xi_i}{i\leq k}$ such that the map $\Phi$ is coded in the model $V(\delta,\xi)$, and $X_\xi:=\Un_{\xi\leq\alpha<\nu}\tilde{M_\alpha}$.
We have $\Phi[M]\sub M_\xi$, and thus $[\Phi(y)<x]\in U$ for all elements $y\in M$ and $x\in X_\xi$. 
It follows that $\Psi[X_\xi^{n+1}\x M]\leU s_{\xi}$.

Fix an ordinal number $\alpha<\xi$.
By the inductive assumption, the set
$(X\cup\Fin)^{n}\x M_\alpha\x M$ is $\Men$,
and thus its continuous image $(X\cup\Fin)^{n}\x \tilde{M_\alpha}\x M$ is $\Men$, too.
The set 
\[
Y:=\bigl((X\cup\Fin)^{n+1}\x M\bigr)\sm\bigl( X_\xi^{n+1}\x M\bigr)
\]
is a union of less than $\card{\xi}$-many sets with the $\Men$ property.
Since $\xi<\addmen$, by Lemma~\ref{lem:BSaddMen}, the set $Y$ is $\Men$.
Let $b'\in\roth$ be a function such that $\Psi[Y]\leinf b'$.
Thus, $\Psi[(X\cup\Fin)^{n+1}\x
M]\leinf\max\{s_\xi,b'\}$.

Let $S:=\sset{s_\xi}{\xi<\nu}$.
Since $\card{S}<\fd$, all finite powers of the set $S$ are $\Men$.
The map $(x,s)\mapsto \sset{n}{x(n)\neq s(n)}$ from the product space $X\x S$ to the set $\roth$, is continuous.
Moreover, it is also surjective: Given $x\in
[\mathbb N]^\infty$, let $\xi$ be such that $x\in V(\delta,\xi)$,
i.e., $x\in M_\xi$. Then $x=\sset{n}{z_\xi(x)(n)\neq s_\xi(n)}$ by
our definition of $z_\xi(n)$. Thus, the product space $X\x S$ is
not $\Men$, because $[\mathbb N]^\infty$ is not
$\Men$ and the $\Men$ property is preserved by continuous images.
\epf

\subsection{Unions and products of $\rothb$ sets with strong properties}
A space is \emph{$\rothb$} if, for each sequence $\eseq{\cU}$ of open covers of that space, there are sets $U_1\in\cU_1,U_2\in\cU_2,\dotsc$ such that the family $\sset{U_n}{n\in\bbN}$ covers the space.

\bthm\label{thm:addproth}
Assume that $\cov(\cM)=\fd$ and the cardinal number $\fd$ is regular. 
\be
\item There are two sets of reals whose all finite powers are $\rothb$, but whose union is not $\Sch$.
\item There are two sets of reals whose all finite powers are $\rothb$, but whose product space is not $\Men$.
\ee
\ethm

In order to prove Theorem~\ref{thm:addproth}, we need the following Lemma.
Let $\kappa$ be an uncountable cardinal number. 
A space $X$ is $\kappa$-concentrated, if there is a countable set $D\sub X$ such that, for any open set $U\sub X$ containing $D$, we have $\card{X\sm U}<\kappa$. 

\blem\label{lem:powerR}
Let $U$ be an ultrafilter such that $\bof(U)\leq\cov(\cM)$, and $X$ be a $U$-scale in $\roth$. 
All finite powers of the set $X\cup\Fin$ are $\rothb$.
\elem

\bpf
Every set of cardinality smaller than $\cov(\cM)$ is Rothberger~\cite[Theorem~4.2]{coc2}, and therefore so is every $\cov(\cM)$-concentrated set.
Since $\bof(U)\leq\cov(\cM)$, the set $X\cup \Fin$ has cardinality smaller than $\cov(\cM)$ or it is $\cov(\cM)$-concentrated~\cite[Proposition~4.3]{ST}. 
Thus, in both cases, the set $X\cup\Fin$ is $\rothb$~\cite[Proposition~4]{reclaw}.
Fix a natural number $n>1$, and assume that the set $(X\cup\Fin)^{n-1}$ is $\rothb$. 
Let $\eseq{\cU}$ be a sequence of open covers of the space $(X\cup\Fin)^{n}$.
For natural numbers $i,j\leq n$, let 
\[
Y^i_j:=
\begin{cases}
X\cup\Fin,& \text{ if }i\neq j,\\
\Fin,&\text{ if }i=j, 
\end{cases}
\] 
and $\pi_i$ be the projection from $(X\cup\Fin)^{n}$ onto the $i$-th coordinate.
A countable union of $\rothb$  sets is $\rothb$.
The set $Y:=\Un_{i\leq n}\prod_{j\leq n}Y^i_j$ is $\rothb$, and thus there are sets $U_1\in\cU_1, U_2\in\cU_2, \dotsc$ such that the set $\Un_kU_k$ covers the set $Y$. 
The set $Y':=(X\cup\Fin)^{n}\sm \Un_kU_k$ is a closed subset of the $\Umen$ set $(X\cup\Fin)^{n}$~\cite[Theorem~5.3(2)]{ST}, and thus it is $\Umen$, too.
For each natural number $i\leq n$, as the projection $\pi_i[Y']$ is a $\Umen$ subset of a $U$-scale $X$, it has cardinality less than $\bof(U)$.
Since $Y'$ is a subset of the set $\prod_{i\leq n}\pi_i[Y']$, a set of cardinality less than $\cov(\cM)$, it is $\rothb$~\cite[Theorem~4.2]{coc2}.
There are sets $U_1'\in\cU_1, U_2'\in\cU_2,\dotsc$ such that $Y'\sub\Un_kU_k'$.
Finally, $(X\cup\Fin)^{n}=\Un_n(U_n\cup U_n')$, and thus the product space $(X\cup\Fin)^{n}$ is $\rothb$~\cite[Theorem~A.1]{book}. 
\epf

\begin{proof}[Proof of Theorem~\ref{thm:addproth}]
(1)
Since  the cardinal number $\fd$ is regular, there is an ultrafilter $U$ with $\bof(U)=\fd$~\cite{canjar}. By Lemma~\ref{lem:NonCoh} and the fact that the near  coherence relation is an equivalence relation, there is an ultrafilter $\tU$ non near coherent to the ultrafilter $U$ with $\bof(\tU)=\fd$.
Since $\cov(\cM)=\fd$ and $\cov(\cM)\leq \fr$~\cite[Theorem~5.19]{blass}, we have $\fd\leq\fr$.
Let $X$ and $Y$ be a $U$-scale, and a $\tU$-scale, respectively, obtained from Lemma~\ref{lem:Uscales}.
We have $\bof(U)=\bof(\tU)=\cov(\cM)$.
By Lemma~\ref{lem:powerR}, all finite powers of the sets $X\cup\Fin$ and $Y\cup\Fin$ are Rothberger.
By Lemma~\ref{lem:Uscales}(1), the union $(X\cup\Fin)\cup(Y\cup\Fin)$ is not Scheepers.

(2) Apply sets from (1) and  Lemma~\ref{lem:Uscales}(2).\epf

\brem
Assuming that $\cov(\cM)=\cof(\cM)$ and the cardinal number $\cov(\cM)$ is regular (a stronger assumption than used in Theorem~\ref{thm:addproth}), there are two \emph{hereditarily} Rothberger sets satisfying the statement of Theorem~\ref{thm:addproth}~\cite{Luzincomb}.
\erem
	
\section{Applications to function spaces} 
Let $X$ be a space, and $\Cp(X)$ be the space of continuous real-valued functions on $X$, with the topology of pointwise convergence.
Properties of the space $X$ can describe local properties of the space $\Cp(X)$, and vice versa. 
E.g., the space $\Cp(X)$ is metrizable (or just first countable) if and only if the space $X$ is countable.
We apply results from the previous sections to products of function spaces.
A space $Y$ has \emph{countable fan tightness}~\cite{cft} if, for each point $y\in Y$ and for every sequence $\eseq{U}$ of subsets of the space $Y$ with $y\in\bigcap_n\overline{U_n}$, there are finite sets $F_1\sub U_1, F_2\sub U_2,\dotsc$ such that $y\in\overline{\bigcup_nF_n}$.
If we request that the above sets $\eseq{F}$ are singletons, then the space $Y$ has \emph{countable strong fan tightness}~\cite{sakai}.
A space is \emph{M-separable}~\cite{BBM} if, for every sequence $\eseq{D}$ of dense subsets of the space there are finite sets $F_1\sub D_1, F_2\sub D_2,\dotsc$ such that the union $\Un_nF_n$ is a dense subset of the space.
If we request that, the above sets $\eseq{D}$ are singletons, then the space is \emph{R-separable}~\cite{BBM}.

\bprp\label{prp:Cp}
\mbox{}

\be
\item Assume that $\fd\leq\fr$ and the cardinal number $\fd$ is regular. 
There are sets of reals $X$, $Y$ such that the spaces $\Cp(X)$, $\Cp(Y)$ have countable fan tightness (are M-separable), but the product space $\Cp(X)\x \Cp(Y)$ does not have countable fan tightness (is not M-separable). 
\item In the Blass--Shelah model,
there are sets of reals $X$, $Y$ such that the spaces $\Cp(X)$, $\Cp(Y)$ have countable fan tightness (are M-separable), but the product space $\Cp(X)\x \Cp(Y)$ does not have countable fan tightness (is not M-separable). 
\item Assume that $\cov(\cM)=\fd$ and  the cardinal number  $\fd$ is regular. 
There are sets of reals $X$, $Y$ such that the spaces $\Cp(X)$, $\Cp(Y)$ have countable strong fan tightness (are R-separable), but the product space $\Cp(X)\x \Cp(Y)$ does not have countable fan tightness (is not M-separable).
\ee
\eprp

To prove Proposition~\ref{prp:Cp}, we need the below result of Scheepers~\cite[Theorem~35, Theorem~13]{Sch}.

\bthm[Scheepers~\cite{Sch}]\label{thm:equiv}
Let $X$ be a set of reals.
\be
\item The space $\Cp(X)$ has countable fan tightness if and only if the space $\Cp(X)$ is M-separable.
\item The space $\Cp(X)$ has countable strong fan tightness if and only if the space $\Cp(X)$ is R-separable.\qed
\ee 
\ethm

\begin{proof}[{Proof of Proposition~\ref{prp:Cp}}]
(1) By Theorem~\ref{thm:addpmen}, there are sets of reals $X$, $Y$ whose all finite powers are $\Men$ but their products space $X\x Y$ is not $\Men$. 
Let $X \sqcup Y$ be the topological sum of these sets.
Since the product space $X\x Y$ is a closed subspace of the product space $(X\sqcup Y)^2$, this product is not $\Men$.
The spaces $\Cp(X)$ and $\Cp(Y)$ have countable fan tightness~\cite[Theorem~3.9]{coc2}. The space $\Cp(X\sqcup Y)$ is homeomorphic to the product space $\Cp(X)\x \Cp(Y)$, and thus it does not have countable strong fan tightness.
For M-separability, apply Theorem~\ref{thm:equiv}(1). 

(2) Proceed as in (1), with the exception that the sets $X$ and $Y$ are as in Theorem~\ref{thm:BSMnonProd}. 

(3) By Theorem~\ref{thm:addproth}(3), there are sets of reals $X$ and $Y$ whose all finite powers are $\rothb$, but their product space $X\x Y$ is not $\Men$.
The spaces $\Cp(X)$, $\Cp(Y)$ have countable strong fan tightness~\cite{sakai}.
As in (1) the product space $\Cp(X)\x \Cp(Y)$ does not have countable fan tightness.
For R-separability and M-separability, apply Theorem~\ref{thm:equiv}. 
\epf

\section{Comments and open problems}

Let $\cov(\cN)$ be the minimal cardinality of a family of Lebesgue null subsets of $\roth$ that covers $\roth$, and $\cof(\cN)$ be the minimal cardinality of a cofinal family of Lebesgue null sets, i.e., every Lebesgue null subset of $\roth$ is contained in a member of the family.
Assuming $\cov(\cN)=\fb=\cof(\cN)$, there is a $\Sch$ set of reals whose square is not $\Men$~\cite[Theorem~43]{cbc}. The $\Sch$ sets constructed here are each $U$-Menger for some ultrafilter $U$, and their finite powers are $\Men$.

\bprb
Assume that $\fd\leq\fr$
\be
\item Is there a $\Sch$ set of reals with a non-$\Men$ square?
\item Is there a $\Sch$ set of reals such that, for each ultrafilter $U$, the set is not $\Umen$?
\ee
\eprb 

By Remark~\ref{rem:MenNonSch}, if the $\Men$ and $\Sch$ properties are different, then the $\Men$ property is not productive.

\bprb
Assume that the Menger and Scheepers properties are equivalent.
Does it follow that the $\Men$ property is productive?
\eprb

By Theorem~\ref{thm:BSMnonProd}, in the Blass-Shelah model the $\Men$ property is not productive.
Thus, it is natural to examine whether the $\Men$ and $\Sch$ properties differ in the Blass--Shelah model.

A subset of $\roth$ is a \emph{semifilter} if it closed under almost supersets.
\emph{Semifilter trichotomy} is the statement: for each semifilter $S$, there is a function $h\in\roth$ such that the set $S/h$ is either the filter of cofinite subsets of $\bbN$, or an ultrafilter, or the full semifilter $\roth$.
Semifilter trichotomy is equivalent to the statement $\fu<\fg$.
Semifilter trichotomy implies that the $\Men$ and $\Sch$ properties are equivalent.~\cite[Theorem~3.7]{semtrich}.
In the Miller model (where semifilter trichotomy holds), the $\Men$ property is productive.
These results motivate the following problem.
\bprb
Assume semifilter trichotomy.
Is the $\Men$ property productive?
\eprb

By Theorems~\ref{thm:MNonSch}(2) and~\ref{thm:add<d}, if $\fd\leq \fr$ or $\add(\Men)<\fd$, then the Menger property is not productive.
\bprb
Is the assumption $\fd\leq\fr$ or $\add(\Men)<\fd$ necessary for the the Menger property being not productive?
\eprb

\subsection{Products of $\Men$ topological groups}

\bprp
Assume that the cardinal number $\fd$ is regular.
\be
\item If $\fd\leq\fr$, then there are two topological groups whose all finite powers are $\Men$ but their product space is not $\Men$. 
\item If $\cov(\cM)=\fd$, then there are two topological groups whose all finite powers are $\rothb$ but their product space is not $\Men$.
\ee
\eprp

\bpf
(1) There is a homeomorphic copy $C\sub\PN$ of the Cantor space, that is linearly independent~\cite{vN}.
By Lemma~\ref{lem:Uscales}, there are ultrafilters $U$, $\tilde{U}$, and subsets $X,Y$ of the set $C$ that are productively $\Umen$, and productively $\tUmen$, respectively, whose product space $X\x Y$ is not $\Men$.
For each natural number $n$, the set \[
\tilde{X}_n:=\sset{x_1\oplus\cdots\oplus x_n}{x_1,\dotsc,x_n\in X},
\] a continuous image of a productively $\Umen$ set $X^n$, is productively $\Umen$. 
Thus the group $\tilde{X}=\Un_n \tilde{X}_n$ is productively $\Umen$, too.
Analogously, construct a $Y$ that is productively $\tUmen$.
Thus, all finite powers of the groups $X$ and $Y$ are $\Men$.
Since the set $C$ is linearly independent, the product space $\tilde{X}\x \tilde{Y}=(X\cap C)\x(Y\cap C)$ is a closed subset of the product space $X\x Y $.
Since the product $X\x Y$ is not $\Men$, the product space $\tilde{X}\x \tilde{Y}$ is not $\Men$, too.

(2) Proceed as in~(1).
By Lemma~\ref{lem:powerR}, all finite powers of the groups $X$, $Y$ are $\rothb$.
\epf

\bprb
Is there, consistently, a $\Men$ topological group whose square is not $\Men$?
\eprb

\subsection*{Acknowledgments}

The first and second named authors thank the third named author for his kind hospitality at the Kurt G\"odel Research Center in fall 2016.
We also thank the Center's Director, researchers and staff for the excellent academic and friendly atmosphere.
The third author would like to thank the Austrian Science Fund FWF (Grants I 2374-N35 and
I 3709-N35) for generous support for this research. We are also grateful to the anonymous referee for careful reading of the manuscript.

\end{document}